\documentclass[letterpaper,onecolumn]{ieeeconf}

\IEEEoverridecommandlockouts                              

\usepackage{graphics} 
\usepackage{epsfig} 
\usepackage{mathptmx} 
\usepackage{times} 
\usepackage{amsmath,amssymb,amsfonts}
\usepackage{algorithmic}
\usepackage{graphicx}
\usepackage{textcomp}
\usepackage{xcolor}
\usepackage{booktabs}
\usepackage{amsfonts}
\usepackage{hyperref}
\usepackage{bm}
\usepackage[english]{babel}

\usepackage{amsthm}
\usepackage{mathtools} 
\usepackage{arydshln}
\newtheorem{definition}{Definition}
\newtheorem{theorem}{Theorem}
\newtheorem{lemma}{Lemma}

\theoremstyle{remark}

\newcommand{\R}{\mathbb{R}}

\newcommand{\Z}{\mathbb{Z}}
\newcommand{\C}{\mathbb{C}}
\newcommand{\bmtx}{\begin{bmatrix}}
\newcommand{\emtx}{\end{bmatrix}}
\newcommand{\bsmtx}{\left[ \begin{smallmatrix}} 
\newcommand{\esmtx}{\end{smallmatrix} \right]}
\newcommand{\bmatarray}[1]{\left[\begin{array}{#1}}
\newcommand{\ematarray}{\end{array}\right]}

\title{\LARGE \bf
Robust Regret Control with Uncertainty-Dependent Baseline
}

\author{Jietian Liu$^{1}$ and Peter Seiler$^{1}$
\thanks{$^{1}$Jietian Liu and Peter Seiler are with the Department of Electrical Engineering and Computer Science, University of Michigan,
Emails: {\tt\small jietian@umich.edu} and {\tt\small pseiler@umich.edu}}%
}

\begin{document}

\maketitle
\thispagestyle{empty}
\pagestyle{empty}

\begin{abstract}

This paper proposes a robust regret control framework in which the performance baseline adapts to the realization of system uncertainty. The plant is modeled as a discrete-time, uncertain linear time-invariant system with real-parametric uncertainty. The performance baseline is the optimal non-causal controller constructed with full knowledge of the disturbance and the specific realization of the uncertain plant. We show that a controller achieves robust additive regret relative to this baseline if and only if it satisfies a related, robust $H_\infty$ performance condition on a modified plant.  One technical issue is that the modified plant can, in general, have a complicated nonlinear dependence on the uncertainty.  We use a linear approximation step so that the robust additive regret condition can be recast as a standard $\mu$-synthesis problem. A numerical example is used to demonstrate  the proposed approach.


\end{abstract}

\section{INTRODUCTION}
This paper considers control design for a discrete-time, linear time-invariant (LTI) plant with parametric uncertainty. The goal is to design an output-feedback controller to: (i) robustly stabilize the uncertain plant, and (ii) ensure a quadratic performance cost remains small in the presence of disturbances and uncertainty. Classical optimal control formulations such as 
the $H_2$ and $H_\infty$ norm~\cite{zhou96,dullerud99} minimize the closed-loop gain from disturbance to error with respect to a fixed nominal plant model. This can lead to robustness issues.


An alternative is regret-optimal control, which measures the performance of a causal controller relative to a baseline controller with greater information or capability. Regret is defined as the performance difference between the
closed-loop cost of the causal controller and that of the baseline
controller. Typical measures for this performance difference include
additive regret, multiplicative regret, or competitive ratio. Prior works have considered baselines such as the optimal noncausal controller with full disturbance preview \cite{goel22TAC,goel19PMLR,goel22CDC,goel21PMLR,sabag21ACC,sabag22CDC,goel23TAC,goel23PMLR,goel22arXivGH}, or the best static state-feedback law \cite{agarwal19,hazan16}. Distributionally robust regret-optimal (DR-RO) control \cite{kargin24pmlr,kargin24pmlrIH,hajar23Allerton,taha23CDC} further extends the idea by evaluating regret under worst-case disturbance distributions in a Wasserstein ball. These formulations typically neglect dynamic and parametric model uncertainty. The control design is performed using a known (nominal) plant model.

In this work, we introduce robust regret-optimal control with an uncertainty-dependent baseline. This builds on our prior work in \cite{liu2023robust_arXiv,liu2024robust} which incorporated model uncertainty in the design but used a baseline controller designed on the nominal model. In this paper, we use a baseline controller that is uncertainty-dependent. This baseline has full knowledge of the disturbance sequence and the specific realization of the uncertainty. The regret is then defined relative to the cost achieved by the optimal baseline in this class for each uncertainty realization. This formulation reduces conservatism by comparing against a more realistic benchmark, while still enforcing robustness over the entire uncertainty set.

The main contributions of this paper are as follows: A) We formulate robust regret with an uncertainty-dependent baseline, and  B) We introduce a linearization step to approximately solve the problem formulated in (A). The remainder of the paper has the following structure.  First, the notation and problem formulation are provided in Section~\ref{sec:prelim}. Next, the main results, including the two key contributions specified above, are given in Section~\ref{sec:MainResult}. Section~\ref{sec:Example} provides a simple scalar example to illustrate the approach. Finally, conclusions and future work are summarized in Section~\ref{sec:conc}.


\section{Preliminaries}
\label{sec:prelim}

\subsection{Notation}
\label{sec:Notation}

This subsection reviews basic notation regarding vectors, matrices,
signals and systems.  This material can be found in most standard
texts on signals and systems, e.g. \cite{zhou96,dullerud99}.

Let $\R^n$ and $\R^{n\times m}$ denote the sets of real $n\times 1$
vectors and $n\times m$ matrices, respectively. Similarly, $\C^n$ and
$\C^{n\times m}$ denote the sets of complex vectors and matrices of
given dimensions. The superscripts $\top$ and $*$ denote transpose and
complex conjugate (Hermitian) transpose of a matrix.  Moreover, if
$M\in \C^{n\times n}$ then $M^{-\top}$ denotes
$(M^\top)^{-1}=(M^{-1})^\top$.  The Euclidean norm (2-norm) for a
vector $v\in \C^n$ is defined to be
$\| v\|_2:= \sqrt{ v^* v} = \sqrt{ \sum_{i=1}^n |v_i|^2}$.  The
induced 2-norm for a matrix $M \in \C^{n\times m}$ is defined to be:
\begin{align*} 
  \|M\|_{2\to 2}:=\max_{0\ne d\in \C^m} \frac{\|Md\|_2}{\|d\|_2}.
\end{align*}
The induced 2-norm for a matrix $M$ is equal to the maximum singular
value, i.e. $\|M\|_{2\to 2} = \bar{\sigma}(M)$ (Section 2.8 of
\cite{zhou96}).

The set of integers is denoted by 
$\Z$.  Let $v:\Z \to \R^n$ and $w:\Z \to \R^n$ be real, vector-valued sequences. Note that we will mainly use two-sided sequences defined
from $t=-\infty$ to $t=+\infty$.  Define the inner product
$\langle v,w \rangle : = \sum_{t=-\infty}^\infty v_t^\top w_t$.  The
set $\ell_2$ is an inner product space with sequences $v$ that satisfy
$\langle v,v\rangle < \infty$. The corresponding norm is
$\|v\|_2 := \sqrt{ \langle v,v \rangle}$. Finally, define the
truncation operator $P_T$ as a mapping from a sequence $v$
to another sequence $w=P_T v$ defined by:
\begin{align*}
  w_t:=\left\{
    \begin{array}{ll}
     v_t  & \mbox{if } t \le T \\ 0 & \mbox{otherwise.}
    \end{array} \right.
\end{align*}


Next, consider a discrete-time, LTI system $G$ with the following
state-space model:
\begin{align}
  \label{eq:Gnotation}
  \begin{split}
  x_{t+1} & = A \, x_t + B\, d_t \\
    e_t & = C \, x_t + D\, d_t,
  \end{split}
\end{align}
where $x_t \in \R^{n_x}$ is the state, $d_t \in \R^{n_d}$ is the
input, and $e_t \in \R^{n_e}$ is the output.  A system $G$ is said to
be causal if $P_T G d = P_T G (P_T d)$, i.e. the output up to time $T$
only depends on the input up to time $T$.  The system is said to be
non-causal if it is not causal, i.e. the output can possibly depend on future
values of the input.  

The matrix $A\in \R^{n\times n}$ is said to be Schur stable if the
spectral radius is $<1$.  If $A$ is Schur stable then $G$ is a causal,
stable system.  Hence $G$ maps an input $d\in\ell_2$ to output
$e\in \ell_2$ starting from the initial condition $x_{-\infty}=0$.
The induced $\ell_2$-norm for a stable system $G$ is defined to be:
\begin{align*}
  \|G\|_{2\to 2}:=\max_{0\ne d\in \ell_2} \frac{\|G d\|_2}{\|d\|_2}.
\end{align*}

Finally, the transfer function for \eqref{eq:Gnotation}
is $G(z) = C (z I_{n_x} - A)^{-1} B + D$. The $H_\infty$ norm
for a stable system $G$ is:
\begin{align*}
  \|G\|_\infty :=\max_{\theta \in [0,2\pi] }  
   \bar{\sigma}\left( G(e^{j\theta}) \right).
\end{align*}
The induced $\ell_2$-norm for a stable system $G$ is equal to the
$H_\infty$ norm, i.e. $\|G\|_{2\to 2} = \|G\|_\infty$ (Theorem 2.3.2 of
\cite{dahleh94}).

\subsection{Problem Formulation}
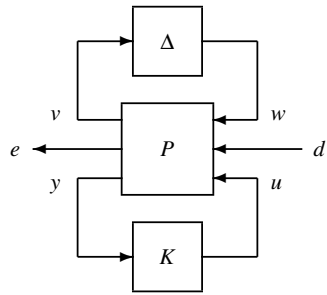
\begin{figure}[h]
\centering
\scalebox{0.85}{
\begin{picture}(140,146)(23,22)
 \thicklines
 \put(75,65){\framebox(40,40){$P$}}
 \put(160,82){$d$}
 \put(155,85){\vector(-1,0){40}}  
 \put(25,82){$e$}
 \put(75,85){\vector(-1,0){40}}  
 \put(80,22){\framebox(30,30){$K$}}
 \put(43,67){$y$}
 \put(55,72){\line(1,0){20}}  
 \put(55,72){\line(0,-1){35}}  
 \put(55,37){\vector(1,0){25}}  
 \put(141,67){$u$}
 \put(135,37){\line(-1,0){25}}  
 \put(135,37){\line(0,1){35}}  
 \put(135,72){\vector(-1,0){20}}  
 \put(80,118){\framebox(30,30){$\Delta$}}
 \put(43,98){$v$}
 \put(55,98){\line(1,0){20}}  
 \put(55,98){\line(0,1){35}}  
 \put(55,133){\vector(1,0){25}}  
 \put(141,98){$w$}
 \put(135,133){\line(-1,0){25}}  
 \put(135,133){\line(0,-1){35}}  
 \put(135,98){\vector(-1,0){20}}  
\end{picture}
} 
\caption{Feedback interconnection $CL(P,K,\Delta)$ for robust synthesis.} 
\label{fig:FLPKUnc}
\end{figure}

Consider the feedback interconnection shown in
Figure~\ref{fig:FLPKUnc}. This interconnection, denoted by $CL(P, K, \Delta)$, consists of a controller $K$ and an uncertainty block $\Delta$ connected to the lower and upper channels of the plant $P$, respectively.  This configuration is standard in robust control and aligns with formulations in Chapter 11 of \cite{zhou96} and Chapter 8 of \cite{dullerud99}. The plant $P$ is a discrete-time, linear time-invariant (LTI) system with additional input/output channels
$(w,v)$ to incorporate the effect of the uncertainty:
\begin{align}
  \label{eq:Punc}
   \bmtx x_{t+1}\\ v_t \\ e_t \\ y_t\emtx =
   \bmtx A & B_w & B_d & B_u\\
   C_v & D_{vw} & D_{vd} & D_{vu} \\
   C_e & 0 & 0 & D_{eu} \\
   C_y & 0 & D_{yd} & 0  \\
   \emtx \bmtx x_t \\ w_t \\ d_t \\u_t \emtx,
\end{align}
where $x_t \in \R^{n_x}$ , $d_t \in \R^{n_d}$ , $e_t \in \R^{n_e}$, $u_t \in \R^{n_u}$, $y_t \in \R^{n_y}$ are the state, disturbance, error, control input, and measured output at time $t$, respectively. $w_t \in \R^{n}$ and $v_t \in \R^{n}$ are the input/output
signals associated with the uncertainty. The
plant~\eqref{eq:Punc} assumes zero feedthrough matrices in several channels, e.g. $D_{yu}=0$ corresponds to no feedthrough from $u$ to $y$.  A standard
loop-shift transformation can be used, under mild technical
conditions, to convert plants with $D_{yu}\ne 0$ and/or $D_{ed} \ne 0$
into the form of Equation~\eqref{eq:Punc} (Section 17.2 of
\cite{zhou96}). We also assume $D_{ew}$ and $D_{yw}$ are zero to simplify the presentation.


Let $\bm{\Delta}$ denote a set of unit-norm bounded, block-diagonal uncertainty matrices in $\R^{n \times n}$ defined by:
\begin{align}\label{eq:UncertaintySet}
\begin{split}
\bm{\Delta} = & \left\{
\mathrm{diag}\left( \delta_1 I_{r_1}, \dots, \delta_S I_{r_S} \right)
\,\middle|\, \right. \left.
\delta_i \in \R,\; |\delta_i| \leq 1 \quad \text{for } i = 1,\dots,S \right\},
\end{split}
\end{align}
where $r_1,\dots,r_s$ are positive integers satisfying $\sum_{i=1}^Sr_i=n$.
This is structured, real parametric uncertainty \cite{zhou96}. The notation $F_U(P,\Delta)$ refers to the system obtained by closing $\Delta \in \bm{\Delta}$ around the top channels of P\footnote{The system $F_L( F_U(P,\Delta), K)$ corresponds to
closing a controller $K$ around the bottom channels of  $F_U(P,\Delta)$.  Thus $F_L( F_U(P,\Delta), K)$ is equivalent to the  shorter notation $CL(P,K,\Delta)$.}.  The set of all possible models is denoted as $\mathcal{M}:=\{ F_U(P,\Delta) \, : \, \Delta \in \bm{\Delta} \}$

The closed-loop system $CL(P,K,\Delta)$ in
Figure~\ref{fig:FLPKUnc} is said to be well-posed if the dynamics have
a unique solution for any disturbance $d \in \ell_2$ starting from
zero initial conditions for $P$ and $\Delta$, i.e. $x_{-\infty}=0$ and
$x^\Delta_{-\infty}=0$. It is common to formulate optimal control problems using one-sided $\ell_2$ signals starting from $t = 0$ with the initial condition $x_0 = 0$. However, a non-causal controller will be introduced
later. Two-sided signals are used to avoid nonzero initial conditions arising from this non-causal controller. 

The cost achieved by a controller $K$ with a disturbance $d \in \ell_2$ and uncertainty
$\Delta\in\bm{\Delta}$ is:
\begin{align}
  \label{eq:JKdDelta}
  J(K,d,\Delta):= \| e\|_2^2 = \sum_{t=-\infty}^\infty e_t^\top e_t.
\end{align}




We consider the performance achieved by an "omniscient" baseline controller that has knowledge of the complete (future) disturbance sequence $d$ and the uncertainty realization $\Delta \in \bm{\Delta}$. This optimal, non-causal baseline controller is denoted $K^{nc}_{d,\Delta}$. This baseline achieves a cost, for the given disturbance $d$ and uncertainty $\Delta$, given by $J(K^{nc}_{d,\Delta},d,\Delta)$. This baseline controller is discussed further in Section~\ref{sec:MainResult}.

We will design our robust controller to minimize the regret relative to this baseline.  The regret is formally defined next.


\begin{definition}
\label{def:robregret}
Let $\gamma_R > 0$ be given. A controller $K$
achieves robust $\gamma_R$-regret relative to the optimal
non-causal controller $K^{nc}_{d,\Delta}$ if 
$CL(P,K,\Delta)$ is well-posed, stable for all
$\Delta \in \bm{\Delta}$, and:
\begin{align}
  \label{eq:robregret}
  \begin{split}
  J(K,d,\Delta)-J(K^{nc}_{d,\Delta},d,\Delta) < \gamma_R^2\, \| d \|_2^2 
  \quad\forall d \in \ell_2, \, d\ne 0 \,\,
  \mbox{ and } \,\, \forall \Delta \in \bm{\Delta}.
  \end{split}
\end{align}
\end{definition}

Definition 1 is a generalization of the more standard (nominal) additive regret. Specifically, the nominal (or best-guess) model is given with $\Delta=0$. The nominal, additive $\gamma_R$-regret feasibility problem is:  Given $\gamma_R>0$, find a controller $K$ such
  that $J(K,d,0)-J(K^{nc}_{d,0},d,0) < \gamma_R^2 \|d\|_2^2$ for all nonzero
  $d\in \ell_2$, or verify that this level of performance cannot be  achieved~\cite{sabag21ACC,goel22CDC}.  
  The $\gamma_R$-regret can be minimized to within any desired tolerance
  using bisection with the $\gamma_R$-regret
  feasibility problem. The $\gamma_R$-regret is additive in the sense
  that $J(K,d,0)$ is within an additive factor $\gamma_R^2 \|d\|_2^2$ of
  the optimal non-causal cost.

The robust regret, as given in Definition 1, uses a baseline controller $K^{nc}_{d,\Delta}$ that has full knowledge of the disturbance and the uncertainty.  This uncertainty-dependent baseline provides the best possible performance for each $\Delta \in \bm{\Delta}$ and each $d\in \ell_2$.  This is in contrast to our prior work on robust regret-optimal control [7] which used $J(K^{nc}_{d,0},d,0)$ as the performance baseline. In this prior work, $K^{nc}_{d,0}$ is the  optimal non-causal controller with full knowledge of $d\in \ell_2$ but designed on the nominal plant.  The baseline used in our current paper, $K^{nc}_{d,\Delta}$, has full knowledge of $d$ and $\Delta$ and hence provides a more meaningful benchmark. This avoids underestimating regret and ensuring proper assessment of robustness.  The next section provides a method to approximately solve the robust regret feasibility problem given in Definition 1.

\section{Main Results}\label{sec:MainResult}
This section provides a method to approximately solve the following
robust $\gamma_R$-regret feasibility problem: Given
$\gamma_R$, find a causal, output-feedback controller $K$
that achieves robust $\gamma_R$-regret relative to $K^{nc}_{d,\Delta}$ or
verify that this level of robust regret cannot be achieved. 

\subsection{Uncertainty-Dependent Non-causal Controller}\label{sec:UncKnc}

We first review the construction of the optimal non-causal controller for the nominal case  ($\Delta=0$). This controller $K^{nc}_{d,0}$ is assumed to have full knowledge of
the plant dynamics, plant state and the (past, current and future)
values of the disturbance.  The controller $K^{nc}_{d,0}$ is optimal in the sense that it
minimizes $J(K,d,0)$ for each $d\in \ell_2$.  A solution for the optimal
non-causal controller is given in Theorem 11.2.1 of
\cite{hassibi99}. The controller is expressed as an operator with
similar results used in \cite{sabag21ACC}.
An explicit state-space model for the finite-horizon, non-causal
controller is constructed in \cite{goel22CDC} using
dynamic programming. The next theorem summarizes these
results and provides a state-space model of $K^{nc}_{d,0}$ for the infinite-horizon case using the stabilizing solution $X$ of a discrete-time algebraic Riccati equation (DARE).  This result follows directly from Lemma 1 and Theorem 1 in \cite{liu2024robust}.

\begin{theorem}
\label{thm:NC}
Let $(A,B_u,C_e,D_{eu})$ be given and define $Q:=C_e^\top C_e$,
$S:=C_e^\top D_{eu}$, and $R:=D_{eu}^\top D_{eu}$.  

Assume: (i) $R \succ 0$, (ii) $(A,B_u)$ is stabilizable, (iii)
$A-B_uR^{-1}S^\top$ is nonsingular, and (iv)
$\bsmtx A-e^{j\theta} I & B_u \\ C_e & D_{eu} \esmtx$ has full column
rank for all $\theta \in [0,2\pi]$. Then:
\begin{enumerate}
\item There is a unique stabilizing
solution $X \succeq 0$ such that $X$ satisfies the following DARE:
\begin{align}
\label{eq:DARE}
  0 & = X  - A^\top XA  - Q  + (A^\top XB_u+S)\, (R+B_u^\top XB_u)^{-1} \, (A^\top XB_u+S)^\top.
\end{align}
\item Let the gain
$K_x:=(R+B_u^\top X B_u)^{-1}(A^\top X B_u+S)^\top$.
Then, $A-B_uK_x$ is a Schur matrix and nonsingular.
\item  Define a non-causal controller $K^{nc}_{d,0}$ with inputs $(x_t,d_t)$ and
output $u_t^{nc}$ by the following update equations:
\begin{align}
  \label{eq:Knc}
  \begin{split}
    v_t & = (A-B_u K_x)^\top ( v_{t+1} + X B_d d_t ), 
    \,\,\, v_\infty=0 \\
    u^{nc}_t & =  -K_x x_t - K_v v_{t+1} - K_d  d_t,
  \end{split}
\end{align}
where 
\begin{align*}
  K_v & := (R+B_u^\top X B_u)^{-1} B_u^\top, \\
  K_d & := (R+B_u^\top X B_u)^{-1} B_u^\top  X B_d.
\end{align*}
Then $J(K^{nc}_{d,0},d,0) \le J(K,d,0)$ for any stabilizing controller $K$ and
disturbance $d\in\ell_2$.
\end{enumerate}
\end{theorem}

Next, consider the uncertain case
with the set of models $\mathcal{M}$. For each $\Delta \in \bm{\Delta}$, the optimal non-causal controller is designed for the uncertain model $F_U(P,\Delta)$, where the specific uncertainty $\Delta$ is known to the controller. The uncertain model $F_U(P,\Delta)$ can be written as:
\begin{align}
  \label{eq:PuncRe}
   \bmtx x^\Delta_{t+1} \\ e^\Delta_t \\ y^\Delta_t \emtx =
   \bmtx A^\Delta & B_d^\Delta & B_u^\Delta \\
         C_e^\Delta & 0 & D_{eu}^\Delta \\
         C_y^\Delta & D_{yd}^\Delta & 0 \emtx
   \bmtx x^\Delta_t \\ d_t \\ u_t \emtx.
\end{align}
Recall that we assumed $P$ has $D_{yw}=0$ and $D_{ew}=0$ in~\eqref{eq:Punc}. This ensures the uncertain model in Equation~\eqref{eq:PuncRe} has $D^{\Delta}_{ed}=0$ and $D^{\Delta}_{yu}=0$ for all $\Delta \in \bm{\Delta}$. Hence we can apply Theorem~\ref{thm:NC} (which assumes this special feedthrough structure) to construct the optimal non-causal  for each plant in the set $\mathcal{M}$. These assumptions on the feedthrough matrices can be relaxed but this requires additional details that will be avoided here.

Let us define: $Q^\Delta := {C_e^\Delta}^\top C_e^\Delta$, $S^\Delta := {C_e^\Delta}^\top D_{eu}^\Delta$, and $R^\Delta := {D_{eu}^\Delta}^\top D_{eu}^\Delta$. Assume that conditions (i)-(iv) in Theorem~\ref{thm:NC} hold for these cost matrices and the uncertain plant $F_U(P,\Delta)$.
Then for each $\Delta \in \bm{\Delta}$, Theorem~\ref{thm:NC} can be applied to obtain the corresponding non-causal controller $K^{nc}_{d,\Delta}$ that minimizes $J(K, d, \Delta)$ for any $d \in \ell_2$.

\subsection{Output Feedback Control Design}
\label{sec:ControlDesign}
We now return to the robust $\gamma_R$-regret feasibility problem:
Given $\gamma_R$, find a causal, output-feedback controller
$K$ that achieves $\gamma_R$-regret relative to $K^{nc}_{d,\Delta}$ or
verify that this level of regret cannot be achieved. This section
provides a conceptual solution to this feasibility problem. Again, we follow the
basic procedure in
\cite{goel22TAC,goel22CDC,goel21PMLR,sabag21ACC}
and use a spectral factorization to reduce the problem to an equivalent robust synthesis problem.

The $\gamma_R$-regret bound in Definition~\ref{def:robregret} 
can be rewritten as:
\begin{align*}
J(K,d,\Delta) & < \gamma_R^2 \, \| d \|_2^2 +  J(K^{nc}_{d,\Delta},d,\Delta) 
 =  \sum_{t=-\infty}^\infty   \gamma_R^2 \, d_t^\top d_t 
  +  (e_{t,\Delta}^{nc})^\top e_{t,\Delta}^{nc}.
\end{align*}
Here $e_\Delta^{nc}$ is the error generated by the closed-loop dynamics with the
non-causal controller $K^{nc}_{d,\Delta}$ after applying Theorem~\ref{thm:NC} for a specific $\Delta \in \bm{\Delta}$. 
Next, define an augmented error as
$\hat{e}_{t,\Delta} := \bsmtx  e_{t,\Delta}^{nc} \\ \gamma_R d_t \esmtx$.  The regret bound in Definition~\ref{def:robregret} can then be written as:
\begin{align}
\label{eq:RegretBnd}
 J(K,d,\Delta) < \sum_{t=-\infty}^\infty \hat{e}_{t,\Delta}^\top \hat{e}_{t,\Delta} = \| \hat{e}_{\Delta}\|_2^2
\end{align}
Details of this process are shown in Section~2.4 in~\cite{liu2024robust}.
Since $K^{nc}_{d,\Delta}$ is non-causal, the right-hand side of \eqref{eq:RegretBnd} corresponds to the squared $\ell_2$ norm of the output of a system mapping $d$ to $\hat{e}_{\Delta}$, which is stable but contains both causal and
non-causal dynamics.  A spectral factorization can be used to re-write
$\| \hat{e}_{\Delta} \|_2^2$ using only stable, causal dynamics.

\begin{lemma}
\label{lem:RegretSF}
Assume that $(A^\Delta,B_u^\Delta,C_e^\Delta,D_{eu}^\Delta)$ satisfy conditions $(i)-(iv)$ in
Lemma~\ref{thm:NC} so that the DARE \eqref{eq:DARE} has a
stabilizing solution $X^\Delta\succeq 0$.

If $\gamma_R>0$ and $((A^{\Delta} - B_u^{\Delta} K_x^{\Delta})^{-\top},X^\Delta B_d^\Delta)$ is stabilizable
then there exists a square $n_d\times n_d$ LTI system $F_\Delta$ such that:
(i) $\| \hat{e}_\Delta\|_2^2 = \| F_\Delta d\|_2^2$ for all $d \in \ell_2$, (ii) $F_\Delta$
is stable, causal, and invertible, and (iii) $F_\Delta^{-1}$ is square,
stable and causal.
\end{lemma}

This corresponds to Lemma 2 in \cite{liu2024robust}.
 It is important to emphasize that the spectral factorization theorem does not
state that $\hat{e}_{\Delta}$ can be computed causally. Rather, it states that
the cost $\|\hat{e}_{\Delta}\|_2^2$ can be equivalently computed as
$\|F_{\Delta} d\|_2^2$ where $F_{\Delta}$ is stable and causal. 
A proof and numerical construction for the spectral factor $F$ can be found in  \cite{liu2023robust_arXiv}.

Using this result, the $\gamma_R$-regret feasibility problem
can be reduced to a robust synthesis problem as follows. The regret in
Definition~\ref{def:robregret} can be written as in
\eqref{eq:RegretBnd}: $J(K,d,\Delta) < \| \hat{e}_{\Delta}\|_2^2$. The left side $J(K,d,\Delta)$ is equal to the closed-loop cost
$\|e\|_2^2 = \|CL(P,K,\Delta) \, d\|_2^2$ achieved with controller $K$
and uncertainty $\Delta \in \bm{\Delta}$.  The bound on the right side is $\|F_\Delta d\|_2^2$ where the spectral factor $F_{\Delta}$ is stable with a
stable inverse.  Define $\hat{d}=F_{\Delta}d$ so that $d=F_{\Delta}^{-1}\hat{d}$.  The
set of all $d\in \ell_2$ maps 1-to-1 on the set of all $\hat{d}\in \ell_2$.
Thus we can rewrite the robust regret bound as
\begin{align}
  \label{eq:robregretF}
  \begin{split}
    \|CL(P,K,\Delta) F_{\Delta}^{-1} \hat{d} \|_2^2  < \| \hat{d} \|_2^2 \quad
    \forall \hat{d} \in \ell_2, \, \hat{d}\ne 0 \,\,
    \mbox{ and } \,\, \forall \Delta \in \bm{\Delta}.
  \end{split}
\end{align}
The system $CL(P,K,\Delta)F^{-1}_{\Delta}$ corresponds to the closed-loop
$CL(P,K,\Delta)$ with the disturbance channel weighted by $F^{-1}_{\Delta}$.
Thus \eqref{eq:robregretF} corresponds to a robust performance
condition on $CL(P,K,\Delta)F^{-1}_{\Delta}$ that must hold for all
uncertainties $\Delta \in \bm{\Delta}$.  In principle, the feasibility of \eqref{eq:robregretF} can be determined using robust synthesis methods, e.g. DK-synthesis~\cite{zhou96,Skogestad05}.
The main technical issue is that $F_\Delta^{-1}$ will have a complicated, nonlinear dependence on $\Delta$. This complication is addressed in the next subsection via a linear approximation.


\subsection{Numerical Approximation}
\label{sec:Numerical}

The robust synthesis problem is difficult to solve since $F^{-1}_{\Delta}$  generally has a nonlinear (and non-rational) dependence on $\Delta$. Therefore, we approximate it using an expansion of the following form:
\begin{align}
    \label{eq:linearize}
    F^{-1}_{\Delta} \approx N_0 + \sum_{i=1}^{S} \delta_i \, N_i.
\end{align}
We construct $\{N_i \}_{i=0}^S$ using a linear fit over the range of uncertainty.
Specifically, let $E_i\in \bm{\Delta}$ denote the uncertainty obtained by setting $\delta_i=1$ and $\delta_j = 0$ for $j\ne i$. Then we define linear expansion with 
$N_0:=F^{-1}_{0}$ and $N_i := \frac{1}{2} (F^{-1}_{E_i}-F^{-1}_{-E_i})$ for $i=1,\ldots, S$.



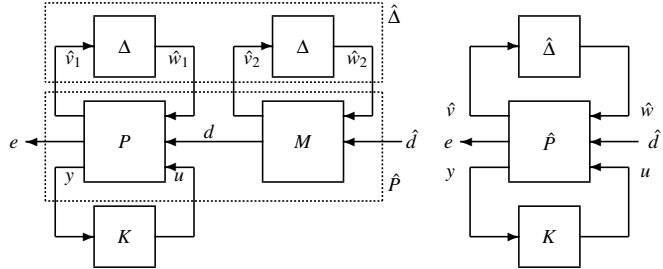
\begin{figure}[h]
\centering
\scalebox{0.75}{
\begin{picture}(200,132)(15,22)
 \thicklines
 \put(50,65){\framebox(40,40){$P$}}
 \put(110,87){$d$}
 \put(140,85){\vector(-1,0){50}}  
 \put(12,82){$e$}
 \put(50,85){\vector(-1,0){30}}  
 \put(55,22){\framebox(30,30){$K$}}
 \put(40,65){$y$}
 \put(35,72){\line(1,0){15}}  
 \put(35,72){\line(0,-1){35}}  
 \put(35,37){\vector(1,0){20}}  
 \put(95,65){$u$}
 \put(105,37){\line(-1,0){20}}  
 \put(105,37){\line(0,1){35}}  
 \put(105,72){\vector(-1,0){15}}  
 \put(55,118){\framebox(30,30){$\Delta$}}
 \put(40,124){$\hat v_1$}
 \put(35,98){\line(1,0){15}}  
 \put(35,98){\line(0,1){35}}  
 \put(35,133){\vector(1,0){20}}  
 \put(92,124){$\hat w_1$}
 \put(105,133){\line(-1,0){20}}  
 \put(105,133){\line(0,-1){35}}  
 \put(105,98){\vector(-1,0){15}}  
 
 \put(140,65){\framebox(40,40){$M$}}
 \put(210,85){\vector(-1,0){30}}  
 \put(212,82){$\hat{d}$}
 \put(145,118){\framebox(30,30){$\Delta$}}
 \put(130,124){$\hat v_2$}
 \put(125,98){\line(1,0){15}}  
 \put(125,98){\line(0,1){35}}  
 \put(125,133){\vector(1,0){20}}  
 \put(182,124){$\hat w_2$}
 \put(195,133){\line(-1,0){20}}  
 \put(195,133){\line(0,-1){35}}  
 \put(195,98){\vector(-1,0){15}} 
 
 \put(30,55){\dashbox(170,55)}
 \put(203,60){$\hat{P}$} 
 \put(30,115){\dashbox(170,40)}
 \put(203,145){$\hat{\Delta}$}
\end{picture}

\hspace{0.1in}

\begin{picture}(110,126)(40,22)
 \thicklines
 \put(75,65){\framebox(40,40){$\hat{P}$}}
 \put(145,82){$\hat{d}$}
 \put(140,85){\vector(-1,0){25}}  
 \put(42,82){$e$}
 \put(75,85){\vector(-1,0){25}}  
 \put(80,22){\framebox(30,30){$K$}}
 \put(43,67){$y$}
 \put(55,72){\line(1,0){20}}  
 \put(55,72){\line(0,-1){35}}  
 \put(55,37){\vector(1,0){25}}  
 \put(141,67){$u$}
 \put(135,37){\line(-1,0){25}}  
 \put(135,37){\line(0,1){35}}  
 \put(135,72){\vector(-1,0){20}}  
 \put(80,118){\framebox(30,30){$\hat\Delta$}}
 \put(43,98){$\hat v$}
 \put(55,98){\line(1,0){20}}  
 \put(55,98){\line(0,1){35}}  
 \put(55,133){\vector(1,0){25}}  
 \put(141,98){$\hat w$}
 \put(135,133){\line(-1,0){25}}  
 \put(135,133){\line(0,-1){35}}  
 \put(135,98){\vector(-1,0){20}}  
\end{picture}
} 
\caption{Approximate representation of $CL(P,K,\Delta)F_\Delta^{-1}$ (Left). This
serial interconnection is in the standard form for robust synthesis with model $\hat{P}$ and uncertainty $\hat{\Delta}$ (Right).} 
\label{fig:SerialConnection}
\end{figure}

An important point is that the expansion in \eqref{eq:linearize} 
can be written as a linear fractional transformation (LFT) $F_U(M,\Delta)$
with $M$ defined as:
\begin{align}
M=
\left[
\begin{array}{ccc: c}
0      & \cdots & 0      & I_{r_1} \\
\vdots & \ddots & \vdots & \vdots \\
0      & \cdots & 0      & I_{r_S} \\
\hdashline
N_1    & \dots  & N_S    & N_0
\end{array}
\right]:=\left[\begin{array}{c: c} M_{11} & M_{12} \\
        \hdashline
         M_{21} & M_{22} 
         \end{array}\right]
\end{align}
The system $CL(P,K,\Delta) F_\Delta^{-1}$, which appears in the robust regret bound   \eqref{eq:robregretF}, can then be approximated by the serial interconnection of 
$CL(P,K,\Delta)$ and $F_U(M,\Delta)$  as shown on the left of  Figure~\ref{fig:SerialConnection}. 

This serial interconnection can be expressed in the standard for for robust control $CL(\hat{P},K,\hat{\Delta})$ with model $\hat{P}$ and 
uncertainty $\hat\Delta$ as shown on the right of Figure~\ref{fig:SerialConnection}.  
The augmented uncertainty $\hat \Delta = \bsmtx \Delta & 0 \\ 0 & \Delta \esmtx$  contains two copies of the original uncertainty.  The augmented plant $\hat{P}$ contains both the original plant $P$ and the model $M$ containing the linear approximation of $F_\Delta^{-1}$.
Finally, the robust regret bound is approximated by:
\begin{align}
  \label{eq:robregretFApprox}
  \begin{split}
    \|CL(\hat{P},K, \hat\Delta) \hat{d} \|_2^2  < \| \hat{d} \|_2^2 \quad \forall \hat{d} \in \ell_2, \, \hat{d}\ne 0 \,\,
    \mbox{ and } \,\, \forall \Delta \in \bm{\Delta}.
  \end{split}
\end{align}
This representation enables the use of standard DK-synthesis methods to compute the robust controller $K$ that satisfies this bound. In
other words, designing a controller $K$ which achieves the regret bound in~\eqref{eq:robregretF} for all admissible $\Delta \in \bm{\Delta}$ has been recast as a  standard robust synthesis problem based on the condition in \eqref{eq:robregretFApprox}.
Consequently, the regret bound is enforced implicitly through the robust $H_\infty$ bound on $CL(\hat{P},K, \hat\Delta)$. This allows established $\mu$-synthesis algorithms to be applied without modification.



\section{Example}
\label{sec:Example}
This section presents an example to illustrate the proposed synthesis method for robust regret control with an uncertainty-dependent baseline. We compare this to nominal and robust regret control with a fixed baseline.\footnote{The code
to reproduce the results is available at:
\url{https://github.com/jliu879/Robust_Regret_Optimal_Control}
}


\subsection{Problem Data}

Consider the following uncertain LTI system
\begin{align}
\begin{split}
  \label{eq:ExampleSystem}
     x_{t+1} &= A_{\Delta}x_t+B_dd_t+B_uu_t\\
     e_t &= \bmtx \sqrt{Q} \\ 0 \emtx x_t+\bmtx 0 \\ \sqrt{R} \emtx u_t
\end{split}
\end{align}
where $A_{\Delta}=A_0+\Delta A_1$ with $\Delta \in \R$ and $ |\Delta| \le 1$. We use the following data in our example:
\begin{align*}
    A_0 = 0.5, \, 
    A_1 = 0.9, \,
    B_d = 5, \,
    B_u = 1, \,
    Q = 3, \, \mbox{ and }
    R = 1.
\end{align*}
Three controllers are synthesized for this problem. First, we compute the standard, nominal additive regret controller $K^{ar}$.  The controller $K^{ar}$ is designed to minimize $\gamma_R$ subject to:
\begin{align*}
    J(K^{ar},d,0) - J(K^{nc}_{d,0},d,0) < \gamma_R^2 \|d\|_2^2 \quad \forall\, d \in \ell_2 \setminus \{0\}.    
\end{align*}
This design uses the nominal plant to construct the baseline controller  $K_{d,0}^{nc}$ and is not robust to the uncertainty.  This problem can be solved using $H_\infty$ synthesis methods on an appropriately weighted plant.  This controller achieves a minimal value (within two digits) of $\gamma_R^{ar} = 0.94$.

Next, we design a robust regret controller with nominal baseline. This controller $K^{rob}$ is designed to minimize $\gamma_R$ subject to:
\begin{align*}
    J(K^{rob},d,\Delta) - J(K^{nc}_{d,0},d,0) < \gamma_R^2 \|d\|_2^2 \quad \forall\, d \in \ell_2 \setminus \{0\},  \forall \Delta \in \bm{\Delta}
\end{align*}
This design also uses the nominal plant to construct the baseline controller $K_{d,0}^{nc}$. It differs from the first design in that it is robust to the uncertainty, i.e. the regret bound holds for all $\Delta \in \bm{\Delta}$.  This can be solved using DK-synthesis methods on an appropriately weighted plant. 
This controller achieves a minimal value of $\gamma_R^{rob} = 3.15$.

Finally, we design a robust regret controller with uncertainty-dependent baseline. This controller $K^{robunc}$ is designed to minimize $\gamma_R$ subject to:
\begin{align*}
    J(K^{robunc},d,\Delta) - J(K^{nc}_{d,\Delta},d,\Delta) < \gamma_R^2 \|d\|_2^2 \quad 
    \forall\, d \in \ell_2 \setminus \{0\},  \forall \Delta \in \bm{\Delta}
\end{align*}
This design uses the uncertain plant to construct the baseline controller $K_{d,\Delta}^{nc}$ and is robust to the uncertainty.  This can be approximately solved using the methods described in this paper.  We use the linear approximation in Section~\ref{sec:Numerical} and DK-synthesis on an appropriately weighted plant. This controller achieves a minimal value of $\gamma_R^{robunc} = 3.75$.


\subsection{Performance Analysis}

\begin{figure}[ht]
    \centering
    \includegraphics[width=0.47\textwidth]{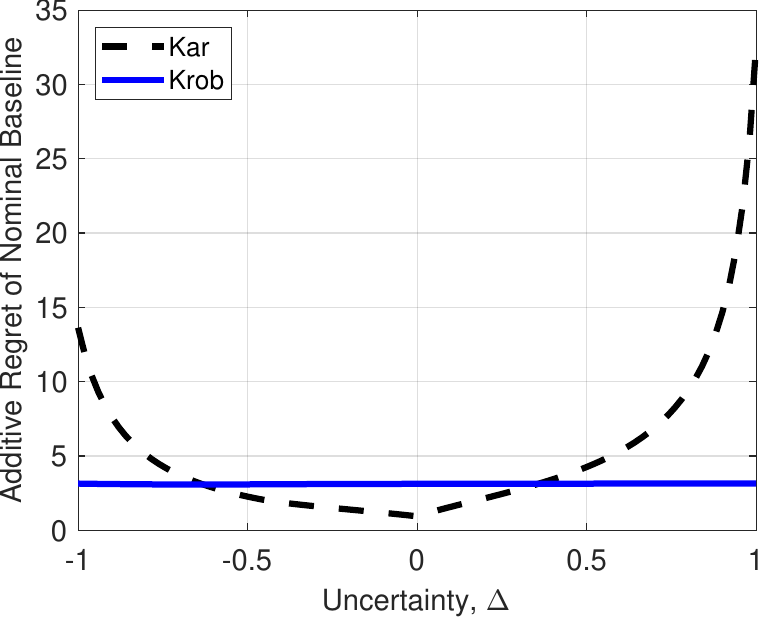}
    \caption{Additive regret using nominal baseline versus uncertainty $\Delta$. }
    \label{fig:regret-nominal}
\end{figure}

We first compare $K^{ar}$ and $K^{rob}$ using the following additive regret:
\begin{align*}
    \text{Regret}(K, \Delta):= \left\{\sup_{\|d\|_2 \leq 1} 
 \left[ J(K,d,\Delta) - J(K^{nc}_{d,0},d,0) \right]\right\}^{1/2}.
\end{align*}
This corresponds to the additive regret measured relative to a fixed baseline $K^{nc}_{d,0}$ designed for the nominal plant, i.e., $\Delta = 0$.\footnote{We compute this regret with the following equivalent expression:
\begin{align*}
\| CL(K,P,\Delta)^* CL(K,P,\Delta) - CL(K^{nc}_{d,0},P,0)^* CL(K^{nc}_{d,0},P,0)\|_\infty^{1/2}.    
\end{align*}}
The performance of $K^{ar}$ and $K^{rob}$ are compared against this fixed baseline for different values of the uncertainty. Fig.~\ref{fig:regret-nominal} shows this additive regret versus uncertainty $\Delta$. The results show that $K_{ar}$ achieves the minimal regret when $\Delta=0$ with a value of 
$\gamma^{ar}_R = 0.94$ (as expected). However, the regret increases significantly as the uncertainty deviates from $\Delta=0$. In contrast, $K_{rob}$ maintains a consistently low regret bounded by
$\gamma_R^{rob} = 3.15$ across the uncertainty range. This confirms its robustness advantage. However,  this formulation may underestimate the true degradation as we use a fixed baseline (with the nominal plant) in this analysis. 

\begin{figure}[ht]
    \centering
    \includegraphics[width=0.47\textwidth]{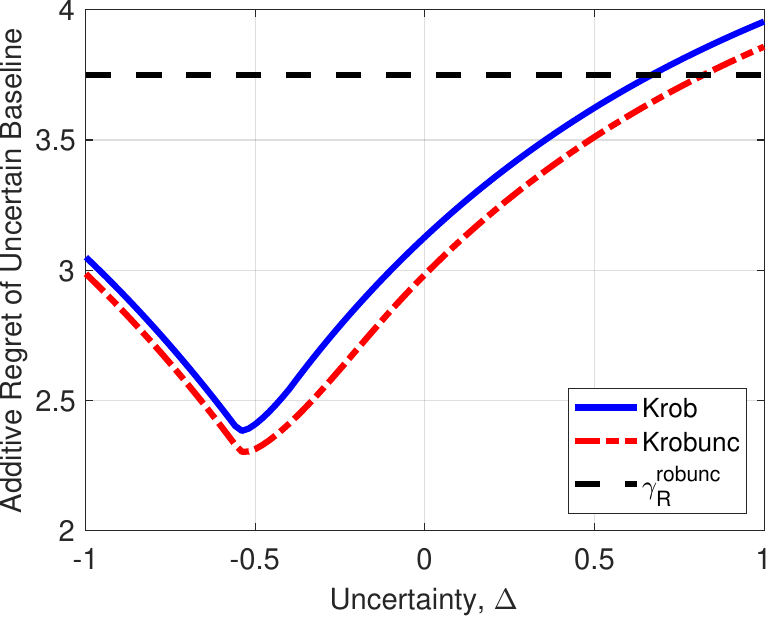}
    \caption{Additive regret using uncertainty-dependent baseline versus uncertainty $\Delta$.}
    \label{fig:regret-uncertain}
\end{figure}

Next, we compare $K^{rob}$ and $K^{robunc}$ using the following additive regret:
\begin{align*}
   \text{Regret}(K, \Delta) := \left\{\sup_{ \|d\|_2 \leq 1} 
 \left[ J(K,d,\Delta) - J(K^{nc}_{d,\Delta},d,\Delta) \right] \right\}^{1/2}.
\end{align*}
This corresponds to the additive regret measured relative to an uncertainty-dependent baseline $K^{nc}_{d,\Delta}$ designed for the uncertain plant. The performance of $K^{rob}$ and $K^{robunc}$ are compared against this uncertainty-dependent baseline for different values of the uncertainty. Fig.~\ref{fig:regret-uncertain} shows this additive regret versus uncertainty $\Delta$. The regret $\gamma^{robunc}_R = 3.75$ is also shown in the plot, which was calculated using linear approximation and DK-synthesis. The  additive regret curve with $K^{robunc}$ (red dash-dot) is expected to lie below the bound $\gamma^{robunc}_R = 3.75$  (black dashed). However, this fails to hold for uncertainty values near 1 due to  the (small) error introduced by the linearized approximation of the spectral factorization.   Improved accuracy can be achieved by employing a higher-order approximation of $F^{-1}_\Delta$.

The uncertainty-dependent baseline formulation better reflects the true performance gap between causal and non-causal control under uncertainty.  In this case, the uncertainty-dependent regret controller ($K^{robunc}$) consistently outperforms the robust regret controller ($K^{rob}$), especially for larger values of $|\Delta|$. This confirms that designing both the controller and the baseline with respect to uncertainty yields more accurate and meaningful regret guarantees. It also highlights the conservative nature of using a fixed nominal baseline in robust settings.

\section{Conclusions}
\label{sec:conc}

This paper introduced a formulation of robust regret control with an uncertainty-dependent baseline. This provides less conservative and more meaningful regret guarantees than fixed nominal baselines.  We also developed a synthesis method that combines spectral factorization, linear approximation, and DK-iteration. Finally, we demonstrated our method on a simple SISO example where the proposed controller outperforms nominal and robust regret controllers under large uncertainty. Future work will focus on replacing the linearization step with a more exact formulations. We will also extend the approach to unstructured or dynamic uncertainty models.





\section*{ACKNOWLEDGMENT}


This material is based upon work supported by the National Science Foundation under Grant No. 2347026. Any opinions, findings, and conclusions or recommendations expressed in this material are those of the author(s) and do not necessarily reflect the views of the National Science Foundation. 

\bibliographystyle{IEEEtran}
\bibliography{reference}

\addtolength{\textheight}{-12cm}   

\end{document}